\documentclass[10pt]{article}
\usepackage[a4paper,
            left=1in,
            right=1in,
            top=1in,
            bottom=1in,
            footskip=.25in
]{geometry}

\usepackage{titlesec}
\usepackage{blindtext,color}
\usepackage{microtype}
\usepackage{scrextend}
\usepackage{upgreek}
\usepackage{stmaryrd}
\usepackage{tocloft}	
\usepackage{multirow}
\usepackage{sansmath}

\usepackage{algorithm}
\usepackage[noend]{algpseudocode}

\usepackage{tcolorbox}
\usepackage{amssymb}
\usepackage{amsmath}
\usepackage{amsthm}
\usepackage{amsfonts}
\usepackage{shadethm}
\usepackage{bbold}
\usepackage{bm}

\usepackage{mathtools}
\DeclarePairedDelimiter\ceil{\lceil}{\rceil}

\usepackage{xspace}
\usepackage{dsfont}
\usepackage{marginnote}
\usepackage{enumerate}% http://ctan.org/pkg/enumerate

\usepackage{graphicx}
\usepackage[citecolor=blue,colorlinks=true,linkcolor=blue]{hyperref}
\usepackage{tikz}
\usetikzlibrary{cd}
\usepackage{cleveref}
\usepackage[colorinlistoftodos,prependcaption,textsize=tiny]{todonotes}

\title{Connecting Kaporin's condition number and the Bregman log determinant divergence}% and low-rank approximations}
\author{Andreas A. Bock \& Martin S. Andersen}
\date{\today}

\newcounter{mythm}
\newcounter{mylem}
\newcounter{myrem}

\newcounter{mypro}
\newcounter{mycol}

\newcounter{mydef}

\newtheorem{theorem}[mythm]{Theorem}
\newtheorem{definition}[mydef]{Definition}
\newtheorem{remark}[myrem]{Remark}
\newtheorem{lemma}[mylem]{Lemma}

\newtheorem{proposition}[mypro]{Proposition}
\newtheorem{corollary}[mycol]{Corollary}

\newcommand{\linearspan}[1]{\textnormal{span}( #1 )}
\newcommand{\Divergence}{\mathcal{D}}
\newcommand{\BregmanLogDet}{\mathcal{D}_\textnormal{LD}}

\newcommand{\cm}{\textsuperscript{\textcolor{red}{[Citation Missing]}}}
\newcommand{\inv}{^{-1}}
\newcommand{\invhalf}{^{-\half}}
\newcommand{\transp}{^\top}
\newcommand{\invtransp}{^{-\top}}
\newcommand{\Hermitian}{^*}
\newcommand{\invHermitian}{^{-*}}

\newcommand{\given}{\,|\,}
\newcommand{\trace}[1]{\Tr\big( #1 \big)}

\newcommand{\logdet}[1]{\log \det( #1 )}
\newcommand{\half}{{\frac{1}{2}}}

\newcommand{\Calpha}{C_\alpha}
\newcommand{\Palpha}{P_\alpha}

\newcommand{\Palphainner}{\Delta_\alpha}
\newcommand{\WV}{\begin{bmatrix} W V \end{bmatrix}}
\newcommand{\WVT}{\begin{bmatrix} W\Hermitian \\ V\Hermitian \end{bmatrix}}

\newcommand{\BCond}[1]{B(#1)}
\newcommand{\KCond}[1]{K(#1)}

\newcommand{\KappaCond}[1]{\kappa_2(#1)}
\newcommand{\IterEpsLTwo}[0]{i_{\ell^2}( \epsilon )}
\newcommand{\IterEpsKaporin}[0]{i_{K}( \epsilon )}

\newcommand{\Sn}{\mathbb{H}^n}

\newcommand{\Snpp}{{\Sn_{++}}}

\newcommand{\sumin}{\sum_{i=1}^n}

\newcommand{\R}{\mathbb{R}}

\newcommand{\Rrr}{\mathbb{R}^{r\times r}}

\newcommand{\Cnr}{\mathbb{C}^{n\times r}}

\newcommand{\lambdaBreg}{\mathfrak{b}}

\DeclareMathOperator{\cone}{cone}

\DeclareMathOperator{\diag}{diag}
\DeclareMathOperator{\Tr}{tr}

\DeclareMathOperator{\rank}{rank}

\DeclareMathOperator*{\subto}{s.t.}

\makeatletter
\newsavebox{\@brx}
\newcommand{\llangle}[1][]{\savebox{\@brx}{\(\m@th{#1\langle}\)}%
  \mathopen{\copy\@brx\mkern2mu\kern-0.9\wd\@brx\usebox{\@brx}}}
\newcommand{\rrangle}[1][]{\savebox{\@brx}{\(\m@th{#1\rangle}\)}%
  \mathclose{\copy\@brx\mkern2mu\kern-0.9\wd\@brx\usebox{\@brx}}}
\makeatother

\newcommand{\BregTruncR}[2]{\llangle #1 \rrangle_{ #2 }}
\newcommand{\BregTrunc}[1]{\BregTruncR{ #1 }{r}}

\algnewcommand\algorithmicoutput{\textbf{Output:}}
\algnewcommand\Output{\item[\algorithmicoutput]}

\DeclareUrlCommand\UScore{}
\newcommand{\expUScore}{%
  \expandafter\expandafter\expandafter
  \UScore
  \expandafter\expandafter\expandafter
}
\begin{document}
\maketitle

\begin{abstract}
\noindent

This paper presents some theoretical results relating the Bregman log
determinant matrix divergence to Kaporin's condition number. These can be viewed
as nearness measures between a preconditioner and a given matrix, and we
show under which conditions these two functions coincide. We also give
examples of constraint sets over which it is equivalent to minimise these two objectives. We focus on preconditioners that are the
sum of a positive definite and low-rank matrix, which were developed in a
previous work. These were constructed as minimisers of the aforementioned divergence, and we show that they are only a constant scaling from also minimising
Kaporin's condition number. We highlight connections to information geometry
and comment on future directions.
\end{abstract}

%\keywords{preconditioner, Kaporin's condition number, Bregman log determinant divergence.}\\

%\tableofcontents

\section{Introduction}

We study preconditioning of the system
\begin{equation}\label{eq:Ax=b}
Ax = b,
\end{equation}
where $A$ is symmetric positive definite.
The purpose of preconditioning is to find a matrix (or operator) $P$ with
inverse $H=P\inv$ so that the matrix
\[
M = H^\half A H^\half
\]
is better conditioned, which can accelerate the progress of iterative methods
such as the conjugate gradient method (CG) \cite{hestenes1952methods}. 
Throughout, we will use the notation above. We also let $I$ denote the identity
matrix, sometimes with a subscript to denote its dimension. 
The \emph{preconditioned} CG method (PCG) for approximating
a solution to $Ax=b$ with a preconditioner $P$ is outlined below:
\begin{subequations}\label{eq:pcg}
\begin{align}
& r_0 = b - Ax_0\\
& p_0 = H r_0\\
& \textnormal{for $k=1, \ldots $}\nonumber\\
& \quad \alpha_k = \frac{r_k\Hermitian H r_k}{p_k\Hermitian A p_k}\\% = \frac{r_k\Hermitian H r_k}{r_{k-1}\Hermitian H\Hermitian A H r_{k-1}}\\
& \quad x_{k+1} = x_k + \alpha_k p_k\label{eq:pcg:xk}\\
& \quad r_{k+1} = r_k - \alpha_k A p_k\\
& \quad \beta_k = \frac{r_{k+1}\Hermitian H r_{k+1}}{r_k\Hermitian H r_k}\\
& \quad p_{k+1} = H r_{k+1} + \beta_k p_k
\end{align}
\end{subequations}
Common termination criteria include an upper limit on the number of iterations,
or the relative residual reaching some prescribed tolerance.\\

\iffalse
It will be useful to recall
\[
\| r_k \|_{A\inv}^2 = \| e_k \|_{A} ^2 = e_k\Hermitian A e_k,
\]
and
\[
\| r_k \|_H^2 = \| e_k \|_{A\Hermitian H A}^2.
\]
\fi
Recently, the authors studied preconditioners obtained as minimisers
of the Bregman
log determinant divergence \cite{bock2023preconditioner,bock2023new}:
\begin{align*}
\BregmanLogDet(A, P) = \trace{AP\inv} - \logdet{AP\inv} - n,
\end{align*}
where feasible preconditioners were selected as a sum of a positive
definite matrix $QQ\Hermitian$ and a low-rank term. We describe this in
greater detail in \cref{sec:precond_low-rank}. The
results in these papers showed that the divergence is capturing a sense of 
nearness to $A$ that may be relevant when searching for a preconditioner $P$
for the PCG algorithm, but without any firm theoretical basis connecting
$\BregmanLogDet(A, P)$ with \cref{eq:pcg}.
In this note, we establish a relationship between this particular divergence
and superlinear convergence of PCG. We highlight the connection between
$\BregmanLogDet(A, P)$ and Kaporin's condition number 
\cite{kaporin1990alternative}.\\

%\subsection{Structure}
This manuscript is structured as follows.
\Cref{sec:spectral} begins by discussing different condition numbers
and relevant literature. These are generally used to characterise the 
progress of \cref{eq:pcg}, which is covered in \cref{sec:convergence_results}.
Next, \cref{eq:Kcond_logdet} contains our main theoretical contributions
related to the condition numbers. As mentioned above, we establish a 
condition on the trace of a preconditioned matrix which implies that
the Bregman log determinant divergence coincides with Kaporin's condition
number. \Cref{sec:precond_low-rank,sec:4.1} 
summarise the low-rank approximation developed in \cite{bock2023new}. In 
\cref{sec:alphacond} we show how the preconditioners associated with these
approximations can be modified to minimise Kaporin's condition number.
\Cref{sec:general_trace_scaling} contains a general result on minimisation
of the Bregman log determinant divergence and Kaporin's condition number.
\Cref{sec:summary} concludes this paper.

%\subsection{Notation}

\section{Spectral and other condition numbers}\label{sec:spectral}

Let $\sigma(A) = \{\sigma_1(A),\ldots,\sigma_n(A)\}$ denote the 
similarly ordered singular values of an $n\times n$ matrix $A$. When
$A$ has real eigenvalues, we denote these by $\lambda(A) = \{\lambda_1(A),\ldots,\lambda_n(A)\}$, with $\lambda_1(A) \geq
\ldots \geq \lambda_n(A)$.\\
%In this section, we describe several condition numbers, their key properties
%and some relations between them.\\

The condition number of a symmetric positive definite matrix $M$ is defined as
\begin{equation}\label{eq:cond:l2}
\KappaCond{M} = \frac{\sigma_1(M)}{\sigma_n(M)}= \frac{\lambda_1(M)}{\lambda_n(M)}.
\end{equation}
This quantity is sometimes used as a basis for discussion of what is a good preconditioner,
when $M=P\invhalf A P\invhalf$, i.e. a matrix $A$ symmetrically preconditioned
by $P$.
% $B$ and $K$
An alternative quantity originated with Kaporin in a series of papers 
\cite{kaporin1990preconditioning,kaporin1990alternative,kaporin1994new,kaporin1994optimization},
with similar approaches dating back to the 80s \cite{gadjokov1985quasi}. Kaporin
introduced the following quantity as a function of a symmetric positive definite 
matrix $M$:
\begin{equation}\label{eq:Kaporin:B}
\BCond{M} = \frac{\frac1n\trace{M}}{\det(M)^{\frac1n}} = \frac{\frac1n \sum_{i=1}^n\lambda_i(M)}{(\Pi_{i=1}^n \lambda_i(M))^{\frac1n }},
\end{equation}
and is the arithmetic mean of the eigenvalues of $M$ divided by the geometric
mean of the eigenvalues. See \cite{kaporin1992explicitly} for the extension to asymmetric matrices.
We shall refer to \cref{eq:Kaporin:B} as \emph{Kaporin's
function}. It decreases when eigenvalues of $M$ are clustered
and depends only modestly on the smallest eigenvalue if it is well-isolated.
The latter is a property not shared by $\KappaCond{M}$ in \cref{eq:cond:l2}. Kaporin's functional is also
closely related to factored sparse approximate inverses (FSAI) \cite{huckle2003factorized,janna2015fsaipack,kolotilina1993factorized,kolotilina1995factorized,kolotilina1999factorized,yeremin2000factorized}.
\Cref{eq:Kaporin:B} leads to Kaporin's condition number:
\begin{equation}\label{eq:Kaporin:K}
\KCond{M} = \BCond{M}^n,
\end{equation}
which is quasi-convex in $M$ and satisfies the following properties \cite{kaporin1994new}, \cite[Theorem 13.5]{axelsson1996iterative}:
\begin{subequations}\label{eq:Kaporin:properties}
\begin{align}
& \KCond{M} \geq 1,\label{eq:Kaporin:properties:iscond}\\
& \KCond{cM} = \KCond{M}, \quad \forall c>0,\label{eq:Kaporin:properties:cscale}\\
& \KCond{X\inv M X} = \KCond{M},\quad \forall X\in\Snpp, \label{eq:Kaporin:properties:scale_invariance}\\
& \BCond{M} \leq \KappaCond{M} \leq \big(\KappaCond{M}^\half + \KappaCond{M}^{-\half}\big)^2 \leq 4\KCond{M} \label{eq:Kaporin_vs_Kappa}.
\end{align}
\end{subequations}

Here, $\Snpp$ is the positive definite cone.
\Cref{eq:Kaporin:properties:iscond} follows from the arithmetic-geometric mean inequality,
with equality if and only if all the eigenvalues of $M$ are equal.
A similar quantity to $\KCond{M}$ has been studied in the optimisation literature 
\cite{dennis1993sizing,jung2025omega}.
%\Cref{eq:Kaporin:properties:scale_invariance} implies that 
%$\KCond{M} = K(cM)$, $c>0$. 
%Further, \cite{kaporin1990alternative} 
%\andreascomment{We sure of this?}
%presents the  following relations between the spectral and $K$-condition numbers:
%\begin{align*}
%& \KappaCond{M} \leq (\sqrt{\KCond{M}} + \sqrt{\KCond{M} - 1}),\\
%& \KCond{M} \leq (\KappaCond{M} - 1)\ln(\KappaCond{M}),
%\end{align*}
While $\KappaCond{M} = \KappaCond{M\inv}$, in general,
\[
\KCond{M} \neq \KCond{M\inv}.
\]
%It can be verified that \cite[Equation 18]{gustafson1999computational} we have 
%the inequality:
%\begin{equation}
%1 \leq \BCond{M} \leq \KappaCond{M}.
%\end{equation}
We also mention that when $n=2$, then $\BCond{M} = \cos \phi(M)\inv$,
where $\cos \phi(M)$ is the first antieigenvalue of $M$ in the sense of
Gustafson  \cite{gustafson1997operator,gustafson1999computational}. See
\cite{gustafson2003operator} for a connection between operator 
trigonometry and preconditioning.\\

Finally, we introduce the Bregman log determinant divergence between
two matrices $A\in\Snpp$, $P\in\Snpp$ whose eigendecompositions are
given by $A = \mathsf{U}{\Xi}\mathsf{U}\Hermitian$ and $P = \mathsf{V}\Omega\mathsf{V}\Hermitian$
\cite{kulis2009low}:
\begin{align}
\BregmanLogDet(A, P) 
& = \trace{AP\inv} - \logdet{AP\inv} - n\nonumber\\
& = \sum_{i=1}^n \sum_{j=1}^n (\mathsf{u}_i\Hermitian \mathsf{v}_j)^2\Big(\frac{\mathsf{\xi}_i}{\mathsf{\omega}_j} - \ln\big(\frac{\mathsf{\xi}_i}{\mathsf{\omega}_j}\big) - 1 \Big).\label{eq:BregmanLogDetDivergence}
\end{align}
We shall also refer to this as the \emph{Bregman divergence} for
brevity, since \cref{eq:BregmanLogDetDivergence} is the only type 
of Bregman divergence we consider in this paper. See \cite{amari2016information} for
other divergences. While this quantity
is not a condition number, since $\BregmanLogDet(A, P)$ is not always
greater than $1$, we shall use it to express several results related
to the convergence rate of PCG.

\subsection{Some convergence results}\label{sec:convergence_results}

In this section, we present several well-known results related to
PCG in terms of the quantities introduced above. The literature on convergence
rates of PCG is vast, so we cover only the essentials. Much material here is
covered, or expanded upon in great detail, in 
\cite{axelsson1996iterative,greenbaum1997iterative,nevanlinna2012convergence}.\\

Recall the symmetrically preconditioned matrix is denoted
\[
M = P\invhalf A P\invhalf.
\]
% Kappa_2
The arguably most well-known convergence result is in terms of 
\cref{eq:cond:l2}:
\begin{align}\label{eq:convergence:linear}
\frac{\| x - x_k \|_A}{\|x - x_0\|_A}
& \leq \frac2{C^k + C^{-k}} \leq 2 C^k,
\end{align}
where
\[
C=\frac{\sqrt{\KappaCond{M}} -1}{\sqrt{\KappaCond{M}} +1}.
\]
\Cref{eq:convergence:linear} is, in general, a pessimistic estimate.
The associated estimate of PCG iterations $\IterEpsLTwo$ needed for a $\epsilon$-reduction in $\| x - x_0 \|_A$, i.e.,
\[
\| x - x_{\IterEpsLTwo} \|_A \leq \epsilon \| x - x_0 \|_A,
\]
where $x_0$ is an initial guess, is therefore bounded from above by 
\cite[Equation 8]{axelsson2001error}:
\[
\IterEpsLTwo = \ceil{\half \sqrt{\KappaCond{M}} \ln\big(\frac{2}\epsilon\big)}.
\]
It is well-known that CG can exhibit superlinear convergence, and several
successful attempts, both theoretical and qualitative, at describing this
phenomenon have been made since the 50s \cite{hayes19543,van1986rate,beckermann2001superlinear}.
A lesser-known bound can also be
expressed in terms of $\KCond{M}$ \cite[Theorem 3.1]{kaporin1994new}:
\begin{equation}\label{eq:convergence:Kaporin}
%\| x - x_j\|_{A\transp P\inv A} \leq (\KCond{M}^{1/k} - 1)^{k/2} \|  x - x_0\|_{A\transp P\inv A},
\| r_k\|_{P\inv} \leq (\KCond{M}^{1/k} - 1)^{k/2} \| r_0 \|_{P\inv}.
\end{equation}
%where $\| x - x_j\|_{A\transp P\inv A} = \| r_j\|_{P\inv }$ is the preconditioned residual.
\Cref{eq:convergence:Kaporin} is only useful if
\begin{equation}\label{eq:Kaporin:iteration:useful}
\sqrt[n]{\KCond{M}} = \BCond{M} < 2,
\end{equation}
or, equivalently,
\[
\log_2 \KCond{M} < n.
\]
It is shown in \cite[Theorem 3.2]{kaporin1994new} that for any
$k\in\mathbb{N}$ and positive $\beta$ (subject to some mild conditions),
there exists a matrix $P\inv $ with $B(P\invhalf A P\invhalf) = \beta$ and an 
initial residual $r_0$ such that
\begin{equation}\label{eq:Kaporin:unimprove}
\| r_k \|_{P\inv}  = (\beta^{n/k} - 1)^{k/2} \| r_0 \|_{P\inv}.
\end{equation}
It is, in this sense, an \emph{unimprovable} estimate.
%and such a preconditioner $P$ 
%is in this sense \emph{optimal}.
%They also provide a bound on the iteration numbers $k$ of PCG required to satisfy
%\[
%\| r_j \|_H \leq \epsilon \| r_0 \|_H,
%\]
%for some $\epsilon>0$.
The following upper bound on the number of iterations required
for an $\epsilon$-reduction in the residual is given in \cite[Theorem 4.1]{kaporin1994new} and is valid for any $\sigma\geq 2$: 
\begin{equation}\label{eq:Kaporin94:iteration_bound}
\IterEpsKaporin = \ceil{\frac{\sigma \ln(\KCond{M}) + 2\ln(\epsilon\inv)}{\sigma \ln (\sigma) - (\sigma-1) \ln (\sigma-1) }}.
\end{equation}
For $\sigma = 2$ this simplifies to
\[
\IterEpsKaporin \leq \ceil{n\log_2(\BCond{M}) + \log_2(\epsilon\inv)},
\]
which is reportedly quite accurate when $n\log_2(\BCond{M}) \gg \log_2(\epsilon\inv)$,
although it is reported that the choice
\[
\sigma = 2 + \frac{\log(\epsilon\inv)}{\ln(\KCond{M})}
\]
provides a more precise approximation.
As pointed out in \cite{axelsson2001optimizing}, when $\BCond{M} \geq 1 + c$, for some $c>0$, then
\[
\IterEpsKaporin \leq nc + \log_2(\epsilon\inv),
\]
which, in terms of $\BCond{M}$, is a pessimistic a priori estimate on 
$\IterEpsKaporin$.
In \cite{axelsson2000sublinear}, several sublinear and 
superlinear convergence results are reported for a variety of approaches
and condition numbers, one of which is
\begin{equation}\label{eq:convergence:Kaporin:3lnD}
\| x_k - x \|_A \leq \Big(\frac{3 \ln (\KCond{M})}{k}\Big)^{k/2} \| x - x_0 \|_A,  % THM 4.3
\end{equation}
which is asymptotically worse than the bound in \cref{eq:Kaporin:unimprove}.
See \cite{axelsson2001error} for more results on iteration and a posteriori error norm
estimates, and \cite{axelsson1995condition} for the study of the rate of convergence for PCG
based on different condition numbers.
\section{Kaporin's condition number and the log determinant divergence}\label{eq:Kcond_logdet}

We now establish some relations between Kaporin's condition 
number \cref{eq:Kaporin:K} to the Bregman log determinant divergence
\cref{eq:BregmanLogDetDivergence}.

\begin{theorem}\label{thm:Bregman_equals_lnKaporin}
\begin{equation}\label{eq:BregmanLogDetInequality}
\BregmanLogDet(A, P) \geq \ln K(AP\inv), \quad \forall A, P \in\Snpp.
\end{equation}
In addition, when
\begin{equation}\label{eq:P:scaled}
\trace{A P\inv} = n
\end{equation}
we have
\begin{equation}\label{eq:logKaporin_equals_Bregman}
\BregmanLogDet(A, P) = \ln K(P\inv A).
\end{equation}
\end{theorem}
\begin{proof}
The inequality \Cref{eq:BregmanLogDetInequality} follows since
\[
\frac1n \trace{P\inv A} - 1 \geq \ln\Big(\frac1n \trace{P\inv A}\Big).
\]
For equality in \Cref{eq:BregmanLogDetInequality}, we see using \cref{eq:P:scaled}, and $M = P\invhalf A P\invhalf$, that
\begin{align*}
\ln (K(M))\nonumber
%& = \ln (B(M)^n)\\
%& = \ln\Big(\big(\frac1n\frac{\trace{M}}{\det(M)^{1/n}}\big)^n\Big)\\
%& = n\ln\Big(\frac1n\frac{\trace{M}}{\det(M)^{1/n}}\Big)\nonumber\\
& = n\ln(\frac1n\trace{M}) - n\frac1n\ln(\det(M))\\
& = n\ln(\frac1n n) - \ln(\det(M))\nonumber\\
& = - \ln(\det(M))\nonumber\\
& = \BregmanLogDet(A, P).\nonumber
\end{align*}
\end{proof}

The Jacobi scaling of a matrix $A\in\Snpp$ by the diagonal matrix 
$\diag(A)$, $\diag(A)_{ii} = A_{ii}$, $i=1,\ldots,n$, is given by
\[
\diag(A)\invhalf A \diag(A)\invhalf.
\]
This scaling is often used to reduce the condition number of a matrix  \cite{greenbaum1997iterative}. Note that
\[
\trace{\diag(A)\invhalf A \diag(A)\invhalf} = n.
\]
The importance of this scaling in constructing FSAI preconditioners
has also been highlighted \cite[Section 2]{yeremin2000factorized}. One may
therefore view the Bregman divergence as penalising the deviation from the
log Kaporin condition number via the term $\trace{P\inv A} - n$.\\

We can also quantify how close the divergence $\BregmanLogDet(A,P)$ is
to $\ln (K(P\inv A))$, as the following corollary shows.

\begin{corollary}\label{cor:error}
Let $A\in\Snpp$, $P=QQ\Hermitian\in\Snpp$ so 
\begin{subequations}\label{eq:AGE_system}
\begin{align}
A & = Q(I + \tilde E)Q\Hermitian,\label{eq:AGE_system:A}\\
\tilde E & = Q\inv (A - P) Q\invHermitian.\label{eq:AGE_system:Q}
\end{align}
\end{subequations}
Then
\[
\ln (K(P\inv A)) = \BregmanLogDet(A,P) + O\big(\trace{\tilde E}^2\big).
\]
%\[
%\tau :=  = \sum_{i=1}^n \lambda_i(\tilde E) = \Big(\sum_{i=1}^n \lambda_i(Q\inv S %Q\invHermitian) \Big) - n.
%\]
\end{corollary}

\begin{proof}
\begin{align*}
\frac1n \BregmanLogDet(A,P)
%& = \frac1n \Big(\big(\sumin 1 + \lambda_i(\tilde E) - \ln(1 + \lambda_i(\tilde E))\big) - n\Big)\\
%& = \frac1n \Big(\sumin 1 + \lambda_i(\tilde E) - \ln(1 + \lambda_i(\tilde E)) - 1\Big)\\
& = \frac1n \Big(\sumin \lambda_i(\tilde E) - \ln(1 + \lambda_i(\tilde E))\Big)\\
& = \Big(\frac1n \sumin \lambda_i(\tilde E)\Big) - \frac1n \sumin \ln(1 + \lambda_i(\tilde E))\\
& \approx \ln( 1 + \frac1n \sumin \lambda_i(\tilde E)) - \frac1n \ln\det(I + \tilde E)\\
& = \ln(\frac1n \sumin 1 +  \lambda_i(\tilde E)) - \frac1n \ln\det(I + \tilde E)\\
& = \ln(\frac1n \trace{I + \tilde E}) - \frac1n \ln\det(I + \tilde E).
\end{align*}
Therefore
%\begin{align*}
%\BregmanLogDet(A,P) 
%& \approx n \ln(\frac1n \trace{I + \tilde E}) - \ln\det(I + \tilde E)\\
%& = \ln (\KCond{M}),
%\end{align*}
%and
\[
\ln (\KCond{M}) = \BregmanLogDet(A,P) + O(n\inv \trace{\tilde E}^2).
\]
\end{proof}
The approximation used in the proof is of course poor when
$\trace{\tilde E} \gg 0$.\\

%\begin{proof}
%See \cref{app:error_proof}.
%\end{proof}
\iffalse
In exact arithmetic, PCG will converge in a single iteration provided
\begin{equation}\label{eq:M_sigma}
P\inv A = c I, \quad c > 0.
\end{equation}
This motivates the introduction of an arbitrary scaling of a preconditioner $P$ i.e.
\[
c P, \quad c >0.
\]
Using this idea, we obtain a connection between $\BregmanLogDet(A, P)$ and $\KCond{P\inv A}$ in the next section.
\fi

Next, we express certain convergence bounds in terms of the Bregman log determinant 
divergence using \cref{thm:Bregman_equals_lnKaporin}.

\subsection{The log determinant divergence and convergence of PCG}\label{sec:superlinear}

In light of \cref{thm:Bregman_equals_lnKaporin}, we can state some of the
convergence results in \cref{sec:convergence_results} in terms of the 
divergence.

\begin{corollary}\label{cor:Divergence:superlinear}
Let $A\in\Snpp$.
For any $k\in\mathbb{N}$, $1 \leq k \leq n-1$ and $\beta>1$, there exists a
preconditioner $P\in\Snpp$ such that
\begin{align}
& \trace{P\invhalf A P\invhalf} = n,\label{superlinear:trace}\\
& \beta = e^{\frac{\BregmanLogDet(A, P)}n},\nonumber
\end{align}
and an initial residual $r_0$ such that the $k$\textsuperscript{th} residual, $r_k$, satisfies
\begin{align*}
\| r_k \|_{P\inv}
& = (\beta^{n/k} - 1)^{k/2} \| r_0 \|_{P\inv}\\
& = (e^{\frac{\BregmanLogDet(A,P)}k} - 1)^{k/2} \| r_0 \|_{P\inv}.
\end{align*}
\end{corollary}
\begin{proof}
Using \cref{superlinear:trace} we have
\[
\KCond{P\inv A}^{1/k} = \big(e^{\ln \KCond{P\inv A})}\big)^{1/k} = \big(e^{\BregmanLogDet(A,P)}\big)^{1/k} = e^{\frac{\BregmanLogDet(A,P)}k}.
\]
The result follows from \cref{eq:Kaporin:properties:cscale} and \cite[Theorem 3.2]{kaporin1994new}.
\end{proof}

By \cref{eq:Kaporin:iteration:useful}, \cref{cor:Divergence:superlinear}
is nontrivial when
\[
%\log_2 \KCond{HA} < k
%\frac{\ln (\KCond{HA})}{\ln(2)} < k
\BregmanLogDet(A,P) < \ln(2)n.
\]

\Cref{cor:Divergence:superlinear} leads to the following two results, where we
note that the trace condition \cref{superlinear:trace} is not necessary.
\begin{corollary}
Let $A\in\Snpp$ and $P\in\Snpp$ be a preconditioner. The $k$\textsuperscript{th} residual satisfies
\begin{equation}\label{eq:convergence:Divergence}
\| r_k \|_{P\inv} \leq \Big(e^{\frac{\BregmanLogDet(A,P)}k} - 1\Big)^{k/2} \| r_0 \|_{P\inv}.
\end{equation}
\end{corollary}
\begin{proof}
The result follows from \cref{thm:Bregman_equals_lnKaporin} and
\cite[Theorem 3.1]{kaporin1994new}.
\end{proof}

\begin{corollary}
Let $A\in\Snpp$ and $P\in\Snpp$ be a preconditioner.
For $k$ even such that
\[
3 \BregmanLogDet(A, P) \leq k<n
\]
we have
\begin{equation}\label{eq:convergence:Divergence:3lnD}
\| x_k - x \|_{A} \leq \Big(\frac{3 \BregmanLogDet(A, P)}{k}\Big)^{k/2} \| x_0 - x \|_{A}.
\end{equation}
%If $\trace{P\inv A} = n$, then \cref{eq:convergence:Divergence:3lnD} coincides
%with \cref{eq:convergence:Kaporin:3lnD}.
\end{corollary}
\begin{proof}
\Cref{eq:convergence:Divergence:3lnD} follows from \cref{thm:Bregman_equals_lnKaporin}
and \cite[Theorem 4.3]{axelsson2000sublinear}.
\end{proof}

We conclude this section with some bounds on the number of iterations required
to achieve an $\epsilon$-reduction in the initial residual of PCG.

\begin{theorem}
Let $A\in\Snpp$ and $P\in\Snpp$ be a preconditioner.
When $\trace{P\invhalf A P\invhalf} = n$, the number of iterations $i_K(\epsilon)$
needed for an $\epsilon$-reduction in the initial $\|r_0\|_{P\inv}$ satisfies
\begin{align*}
i_K(\epsilon) \leq \ceil{\ln(2)\inv\big(\ln(\epsilon\inv) + D(A,P)\big)}.
\end{align*}
\end{theorem}
\begin{proof}
Using \cref{eq:Kaporin94:iteration_bound} for $\sigma=2$ we have
\begin{align*}
i_K(\epsilon)
& \leq \ceil{n\log_2(\BCond{P\invhalf A P\invhalf}) + \log_2(\epsilon\inv)}\\
& = \ceil{\log_2(\KCond{P\inv A}) + \log_2(\epsilon\inv)}\\
& = \ceil{\log_2(\frac{1}{\det(P\inv A)}) + \log_2(\epsilon\inv)}\\
& = \ceil{\frac{\ln(\frac{1}{\det(P\inv A)})}{\ln(2)} + \frac{\ln(\epsilon\inv)}{\ln(2)}}\\
& = \ceil{\ln(2)\inv\big(\ln(\epsilon\inv) - \ln(\det(P\inv A))\big)}\\
& = \ceil{\ln(2)\inv\big(\ln(\epsilon\inv) + D(A,P)\big)}.
\end{align*}
\end{proof}
%\section{Application to preconditioning}\label{sec:precond_low-rank}
\section{Preconditioners based on low-rank approximations}\label{sec:precond_low-rank}

% Basic intro
%In this section, we discuss preconditioners involving low-rank approximations.

% Kaporin
Kaporin also constructs preconditioners using low-rank matrices
\cite{kaporin1994optimization,kaporin1990preconditioning,kaporin1992two}. In
\cite{kaporin1994optimization}, a sparse lower triangular matrix 
variable $G$ is introduced, and the expression
\[
\BCond{GAG\transp}
\]
is minimised subject to some structural constraints on $G$. This is referred to
as the first stage. Second, they construct a
matrix
\begin{equation}\label{eq:F}
F = I + CSC\transp
\end{equation}
such that $B(F GAG\transp)$ is minimised. Here, some structure is imposed on $C\in\mathbb{R}^{n\times r}$ and the matrix $S\in\mathbb{R}^{r\times r}$ that
minimises $B(F GAG\transp)$ has the following form:
\[
S = \pi_0 C\transp GAG\transp C - (C\transp C)\inv \in\mathbb{R}^{r \times r},
\]
where
\[
\pi_0 = \frac{n - \trace{(C\transp C)\inv C\transp GAG\transp C}}{n - r}.
\]
\cite{axelsson2001optimizing} adopts a similar two-stage approach.
The preconditioner described above can be written as
\begin{equation}\label{eq:Gprec}
G\invtransp F\inv G\inv = G\invtransp (I + CSC\transp)\inv G\inv.
\end{equation}
If we let $G\invtransp G\inv \approx A$ be an incomplete Cholesky factorisation of
$A$, and recall that for some $n\times n$ matrix $W$ with $\rank W = r$,
\[
F\inv = (I + CSC\transp)\inv = I + W,
\]
then \cref{eq:Gprec} can be written as
\[
Q(I + W)Q\transp, \qquad \rank W = r.
\] 
This is similar to the preconditioners sought in
\cite{bock2023preconditioner,bock2023new}. In \cref{eq:F} above, the low-rank
term has a specific structure given by the choice of $C$, which is more restrictive
than the general low-rank constraint considered in the previous references.
In this report, we focus preconditioners for which $\KCond{M}$ or $\BregmanLogDet(A,P)$
is minimised; however, other approaches include $\kappa_2$-based 
minimisation, see e.g. \cite{qu2024optimal,gao2023scalable} and the 
references therein.

\subsection{Bregman log determinant divergence-based preconditioners}\label{sec:4.1}

% Our work
In this section, we describe the preconditioners introduced in
\cite{bock2023preconditioner,bock2023new}, which combine approximate 
factorisations with low-rank approximations sought as minimisers of a 
Bregman divergence. We introduce the approach taken in these 
references to generalise the resulting preconditioners in 
\cref{sec:alphacond}, and establish the connection with
Kaporin's condition number.\\

In \cite{bock2023preconditioner,bock2023new}, the authors assume that an approximate
factorisation $QQ\Hermitian$ of the target matrix $A\in\Snpp$ is assumed available
(e.g. incomplete Cholesky). Writing $A$ in terms of $Q$ we obtain
\begin{equation}\label{eq:A_Qform}
A = Q(I + \tilde E)Q\Hermitian,
\end{equation}
where $\tilde E = Q\inv A Q\invHermitian - I$.
A low-rank term was introduced to produce
the preconditioner
\begin{equation}\label{eq:Bregman_preconditioner_form}
P = Q(I + VDV\Hermitian)Q\Hermitian,
\end{equation}
where $V\in\Cnr$ and $D\in\Rrr$ are chosen such that
\[
VDV\Hermitian \approx \tilde E.% =  Q\inv (A - P) Q\invHermitian = Q\inv A Q\invHermitian - I. 
\]
The approximation above was sought in the sense of the Bregman log
determinant divergence introduced in \cref{eq:BregmanLogDetDivergence}. 
This divergence is invariant to congruence transformations. Indeed,
using \cref{eq:A_Qform,eq:Bregman_preconditioner_form} we obtain
\[
\BregmanLogDet(A, P) = \BregmanLogDet(I + \tilde E, I + VDV\Hermitian).
\]
Then, the matrices $V\in\Cnr$ and $D\in\Rrr$ are found as minimisers of
\begin{subequations}\label{eq:Bregman_preconditioner_pbm}
\begin{align}
\min_{V\in\Cnr,\; D\in\Rrr} \quad & \BregmanLogDet(I + \tilde E, I + VDV\Hermitian)\label{eq:Bregman_preconditioner_pbm:obj}\\
\subto\quad
& I + VDV\Hermitian \in \Snpp\\
& \rank D \leq r.
\end{align}
\end{subequations}
In \cite{bock2023preconditioner}, $\tilde E$ was assumed to be positive 
semidefinite, in which case $V$ and $D$ were found as the constituents
of a truncated singular value decomposition (TSVD) of $\tilde E$.
\cite{bock2023new} generalised the approach and allowed $\tilde E$ to be
indefinite, and it was shown that a low-rank matrix $VDV\Hermitian$ that minimises
\cref{eq:Bregman_preconditioner_pbm:obj} can differ
from a TSVD of $\tilde E$. To see this, let
\begin{equation}\label{eq:E:eigendecomp}
\tilde E = U\Theta U\Hermitian
\end{equation}
be an eigendecomposition of $\tilde E$, and let $\Pi\in\mathbb{R}^{n\times n}$ be a diagonal matrix with positive entries. The divergence $\BregmanLogDet(I + U\Theta U\Hermitian, I + U\Pi U\Hermitian)$ takes the form of the following simple expression
 (cf. \cref{eq:BregmanLogDetDivergence}):
\begin{equation}\label{eq:Divergence:U=V}
\BregmanLogDet(I + U\Theta U\Hermitian, I + U\Pi U\Hermitian) = \sum_{i=1}^n \Big( \frac{1 + \theta_i}{1 + \pi_i} - \log \big(\frac{1 + \theta_i}{1 + \pi_i}\big) - 1\Big)
\end{equation}
since eigenspaces of the inputs are aligned. If $\pi_i = \theta_i$ for any 
$i=1,\ldots,n$, the $i^\text{th}$ summand of \cref{eq:Divergence:U=V} 
vanishes. If we consider \cref{eq:Divergence:U=V} as a function of the 
diagonal matrix $\Pi$ with the constraint $\rank \Pi \leq r$, a minimiser
will satisfy $\pi_i = \theta_i$ for some indices $i$,
and $\pi_i = 0$ otherwise (by the rank constraint). The $n-r$ summands of 
\cref{eq:Divergence:U=V} that do not vanish are therefore of the form
\[
\theta_i - \log (1 + \theta_i).
\]
By seeking a low-rank approximation of $\tilde E$ in the sense of \cref{eq:Bregman_preconditioner_pbm}, it is therefore
of interest to identify the eigenvalues $\theta_i$ that take on the largest value 
under the image of the map
\begin{equation}\label{eq:bregman_curve}
\gamma(\lambda) = \lambda - \log\big(1+\lambda\big).
\end{equation}
$\gamma$ is different from the map $\lambda\mapsto |\lambda|$, which can lead to
$VDV\Hermitian$ being different from a TSVD of $\tilde E$.
%\Cref{eq:bregman_curve} is shown in \cref{fig:bregman_curve}.
%\begin{figure}%[h!]
%    \centering
%    \includegraphics[scale=0.4]{figures/bregman_curve.pdf}
%    \caption{Bregman curve $\gamma$ in \cref{eq:bregman_curve}.}
%    \label{fig:bregman_curve}
%\end{figure}
%Indeed, suppose $U$ is given, and we were to minimise \cref{eq:Divergence:U=V}
%with the restriction that the approximant $\Pi$ has $r<n$ non-zero entries, this places 
%the contribution to \cref{eq:Divergence:U=V} as a result of not choosing a 
%particular eigenvalue $\theta_i$, $i=1,\ldots,n$, is given by $\gamma(\theta_i)$.
The map $\gamma$ induces an order on the eigenvalues of a matrix
$A$, which is captured by the following definition.

\iffalse
\begin{definition}[order induced by $\gamma$]\label{def:order}
Let $\lambdaBreg(A) = \{\lambdaBreg_1(A), \ldots, \lambdaBreg_n(A) \}$ denote
the set of eigenvalues of $A$ with
\begin{equation}\label{eq:partial_order}
\lambdaBreg_1(A) \succeq \ldots \succeq \lambdaBreg_n(A),
\end{equation}
where the order $\succeq$ is given by
\[
\lambdaBreg_i(A) \succeq \lambdaBreg_j(A) \qquad \textnormal{if} \qquad \gamma(\lambdaBreg_i(A)) \geq \gamma(\lambdaBreg_j(A)).
\]
In other words, the index subscript on $\lambdaBreg(A)$ enumerates the eigenvalues
according to their magnitude under $\gamma$ (as opposed to the 
common algebraic ordering, cf. \cref{def:lambda}).
\end{definition}
\Cref{def:order} leads to the definition of \emph{Bregman log determinant truncation}:
\fi

\begin{definition}[Bregman log determinant truncation \cite{bock2023new}]\label{def:Bregman_truncation}
Let $I + X \in \Snpp$ and let
\[
X = W\Lambda W\Hermitian
\]
be an eigendecomposition where the diagonal elements of $\Lambda$ are sorted in algebraically non-increasing order and $W$ is orthonormal. We introduce a permutation matrix $\mathcal{P}$ and define $Z$ and $\mathfrak{B}$ via
\[
X = W \mathcal{P}\transp \mathfrak{B} \mathcal{P} W\Hermitian = Z \mathfrak{B} Z\Hermitian,
\]
where 
\[
\mathfrak{B} = 
\begin{bmatrix}
\lambdaBreg_1(X) & & \\
& \ddots & \\
& & \lambdaBreg_n(X)
\end{bmatrix},
\]
is the diagonal matrix containing the eigenvalues $X$ given by the order $\succeq$:
\[
\lambdaBreg_i(A) \succeq \lambdaBreg_j(A) \qquad \textnormal{if} \qquad \gamma(\lambdaBreg_i(A)) \geq \gamma(\lambdaBreg_j(A)).
\]
Here, $\gamma$ is the function given in \cref{eq:bregman_curve}.
We define a BLD truncation of $X$ to order $r$ by
\begin{equation}\label{eq:Bregman_truncation}
\BregTrunc{X} = \begin{bmatrix} z_1 | \cdots | z_r \end{bmatrix} \begin{bmatrix}
\lambdaBreg_1(X) & & \\
& \ddots & \\
& & \lambdaBreg_r(X)
\end{bmatrix} \begin{bmatrix} z_1 | \cdots | z_r \end{bmatrix}\Hermitian.
\end{equation}
\end{definition}
In other words, the BLD truncation is found by selecting the rows and
columns of an eigendecomposition of a matrix corresponding to the $r$ 
largest values of the eigenvalues under the map $\gamma$.\\
 
\Cref{def:Bregman_truncation} leads to the following preconditioner.
\begin{definition}[Bregman log determinant precondition \cite{bock2023new}]\label{def:Bregman_preconditioner}
Let $A \in \Snpp$, $QQ\Hermitian \in \Snpp$ and $1\leq r < n$, and 
assume $\tilde E$ has rank greater than $r$.
We  call 
\begin{equation}\label{eq:Bregman_preconditioner}
P = Q(I + \BregTrunc{\tilde E})Q\Hermitian.
\end{equation}

a \emph{Bregman log determinant preconditioner}.
\end{definition}

\subsection{Minimisation of the Bregman divergence and Kaporin's condition number}\label{sec:alphacond}

%From what precedes, we know that the preconditioner associated with a
%minimiser of \cref{eq:Bregman_preconditioner_pbm} does not necessarily 
%minimise Kaporin's condition number.
In this section, we modify \cref{eq:Bregman_preconditioner} by permitting
an arbitrary constant scaling of the preconditioner:
\begin{equation}\label{eq:Bregman_preconditioner:alpha}
\Palpha = Q(\alpha(I - VV\Hermitian) + V(I_r + D)V\Hermitian)Q\Hermitian.     
\end{equation}
%We call this a \emph{scaled} Bregman log determinant preconditioner.
When $\alpha=1$, we recover \cref{eq:Bregman_preconditioner}.
We proceed to show that a specific choice of $\alpha$ leads to 
\cref{eq:Bregman_preconditioner:alpha} minimising Kaporin's condition number,
as well as the Bregman divergence. Consider the following expression
\begin{equation}\label{eq:divergence:alpha}
\BregmanLogDet(A, \Palpha) = \BregmanLogDet(I + \tilde E, \alpha(I - VV\Hermitian) + V(I_r + D)V\Hermitian).
\end{equation}
\Cref{thm:alpha} below describes a minimiser of \cref{eq:divergence:alpha} as a function of $\alpha$,
letting $V$ and $D$ be fixed:

\begin{theorem}\label{thm:alpha}
Suppose $V\in\Cnr$, $D\in\Rrr$ have been chosen according to a Bregman truncation i.e.
\[
\BregTrunc{\tilde E} = VDV\Hermitian,
\]
and set
\[
P = Q(I + \BregTrunc{\tilde E})Q\Hermitian.
\]
Then,
\begin{equation}\label{eq:alpha}
\alpha^\star = \frac{\trace{(I + \tilde E)(I - VV\Hermitian)}}{n - r} =
\frac{\trace{P\inv A} - r}{n - r} = \frac{\sum_{i=r+1}^n 1 + \lambdaBreg_i(\tilde E)}{n - r}
\end{equation}
is the unique minimiser of \cref{eq:divergence:alpha} as a function of $\alpha$.
Furthermore,
\begin{equation}\label{eq:BregmanEqualslogKaporin}
\BregmanLogDet(A,\Palpha) = \ln \KCond{\Palpha\inv A} = - \ln \det(\Palpha\inv A) = - \ln\det(P\inv A) + (n-r)\ln(\alpha).
\end{equation}
\end{theorem}
\begin{proof}
Let $\Palphainner = \alpha(I - VV\Hermitian) + V(I_r + D)V\Hermitian$. Then
\[
\frac{\text{d}}{\text{d}\alpha}\BregmanLogDet(A, \Palpha) = - \trace{\Palphainner\inv(I + \tilde E)\Palphainner\inv \frac{\text{d}\Palphainner}{\text{d}\alpha}} + \frac{\text{d}}{\text{d}\alpha} \ln\det \Palphainner.
\]
Since
\[
\frac{\text{d}\Palphainner}{\text{d}\alpha} = I - VV\Hermitian,
\]
\[
- \trace{\Palphainner\inv(I + \tilde E)\Palphainner\inv \frac{\text{d}\Palphainner}{\text{d}\alpha}} = - \alpha^{-2} \trace{(I + \tilde E)(I - VV\Hermitian)}.
\]
Combining this with
\[
\frac{\text{d}}{\text{d}\alpha} \ln\det \Palphainner = \alpha\inv \trace{I - VV\Hermitian}
\]
yields the equation
\[
- \alpha^{-2} \trace{(I + \tilde E)(I - VV\Hermitian)} + \alpha\inv \trace{I - VV\Hermitian} = 0.
\]
\Cref{eq:alpha} follows as a consequence. Uniqueness follows since
$\BregmanLogDet(A, \Palpha)$ is convex in $\beta = \alpha\inv$.\\
By direct computation,
\[
(I + \tilde E)(I - VV\Hermitian)  = \trace{P\inv A} - r.
\]
Finally, \cref{eq:BregmanEqualslogKaporin} follows since
\begin{align*}
\BregmanLogDet(A, \Palpha)
%& = \sum_{i\not\in \BregmanIndexSet{\tilde E}} \frac{1 + \lambda_i(\tilde E)}{\alpha^\star} - \ln\Big(\frac{1 + \lambda_i(\tilde E)}{\alpha^\star} \Big) - 1\\
%& = \sum_{i\not\in \BregmanIndexSet{\tilde E}} \frac{1 + \lambda_i(\tilde E)}{\frac{\trace{P\inv A} - r}{n - r}} - \ln\Big(\frac{1 + \lambda_i(\tilde E)}{\frac{\trace{P\inv A} - r}{n - r}} \Big) - 1\\
& = \frac{\trace{P\inv A} - r}{\frac{\trace{P\inv A} - r}{n - r}} - \ln\det(\Palpha\inv A) - (n - r)\\
& = - \ln\det(\Palpha\inv A).
\end{align*}
\end{proof}

%\begin{definition}[Normalised Bregman log determinant preconditioner]
%We refer to $P_{\alpha^\star}$ (i.e. \cref{eq:Bregman_preconditioner:alpha}
%with $\alpha$ chosen according to \cref{eq:alpha}) as a \emph{normalised} %Bregman
%log determinant preconditioner.
%\end{definition}

$\alpha$ in \cref{eq:Bregman_preconditioner:alpha} is very similar to the
$\sigma$ variable in \cite[Section 3.3]{axelsson2001optimizing}, which was
used to move the set of smallest eigenvalues of the preconditioned matrix
closer to the larger eigenvalues.

%\begin{remark}
%\Cref{eq:BregmanEqualslogKaporin} implies that the bounds \cref{eq:convergence:Kaporin,eq:convergence:Kaporin:3lnD} are
%valid for CG applied to the matrix $\Palpha\inv A$. The difference
%between $P$ and $\Palpha$ is simply the scaling of the remaining
%eigenvalues of $I + \tilde E$ not captured by the low-rank truncation.
%\end{remark}

While the choice \cref{eq:alpha} minimises Kaporin's condition number and
Bregman divergence as a function of $\alpha$, the following theorem
reveals that many different values of $\alpha$ minimise the condition number
i.e. the map
\[
\alpha \mapsto \kappa_2(\Palpha\invhalf A\Palpha\invhalf).
\]

\begin{theorem}\label{thm:alpha:cond}
Suppose $V$ and $D$ have been chosen according to a Bregman truncation
(cf. \cref{def:Bregman_preconditioner}).
Then any 
\begin{equation}\label{eq:alpha:interval} 
\alpha \in [\lambdaBreg_{r+1}(A), \lambdaBreg_n(A)]
%& = [\lambda_1((I + \tilde E)(I-VV\Hermitian)), \lambda_n((I + \tilde E)(I-VV\Hermitian))]
\end{equation}
will be a minimiser of
\[
\alpha \mapsto \kappa_2(\Palpha\invhalf A\Palpha\invhalf).
\]
%and
%\[
%\lambdaBreg_n(P\inv A) \leq 1 \leq \lambdaBreg_{r+1}(P_\alpha\inv A).
%\]
\end{theorem}
\begin{proof}
We have
\begin{equation}
\kappa_2(\Palpha\invhalf A\Palpha\invhalf) = \frac{\max\Big(1, \frac{\lambda_1((I + \tilde E)(I-VV\Hermitian))}{\alpha}\Big)}{\min\Big(1, \frac{\lambda_n((I + \tilde E)(I-VV\Hermitian))}{\alpha}\Big)} \geq \frac{\lambda_1((I + \tilde E)(I-VV\Hermitian))}{\lambda_n((I + \tilde E)(I-VV\Hermitian))} =: \kappa_2^\star(\Palpha\invhalf A\Palpha\invhalf).
\end{equation}
If $\alpha > \lambda_1((I + \tilde E)(I-VV\Hermitian))$,
\[
\kappa_2(\Palpha\invhalf A\Palpha\invhalf) = \frac{\alpha}{\lambda_n((I + \tilde E)(I-VV\Hermitian))} > \kappa_2^\star(\Palpha\invhalf A\Palpha\invhalf).
\]
Conversely, if $\alpha < \lambda_n((I + \tilde E)(I-VV\Hermitian))$,
\[
\kappa_2(\Palpha\invhalf A\Palpha\invhalf) = \frac{\lambda_1((I + \tilde E)(I-VV\Hermitian))}{\alpha} > \kappa_2^\star(\Palpha\invhalf A\Palpha\invhalf).
\]
The result follows by \cref{def:Bregman_truncation}, namely
\[
[\lambdaBreg_{r+1}(A), \lambdaBreg_n(A)] = [\lambda_1((I + \tilde E)(I-VV\Hermitian)), \lambda_n((I + \tilde E)(I-VV\Hermitian))]
\]
\end{proof}

\begin{figure}[h!]
    \centering
    \includegraphics[width=0.5\linewidth]{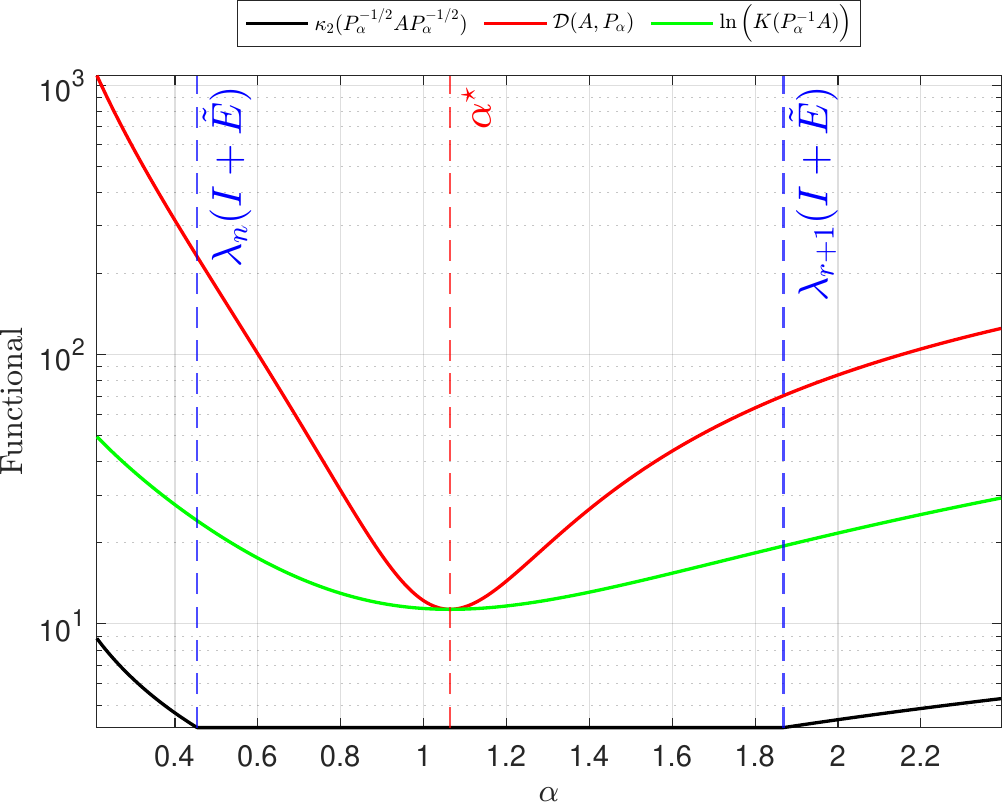}
    \caption{Numerical illustration of how $\alpha$ influences the values $\kappa_2(\Palpha\invhalf A\Palpha\invhalf)$, 
$\BregmanLogDet(A, \Palpha)$, $\ln\big(\KCond{\Palpha\inv A}\big)$.}
    \label{fig:alpha_functionals}
\end{figure}

We conclude this section with a numerical illustration of the effect of $\alpha$
on the condition numbers seen above. $A\in\Snpp$ is the matrix \texttt{494\_bus}
from the SuiteSparse Matrix Collection \cite{kolodziej2019suitesparse} of order
$n=494$ and $Q$ given by a zero fill incomplete Cholesky decomposition of $A$. 
$\Palpha$ is defined in \cref{eq:Bregman_preconditioner:alpha}, with $r=49$.
In \cref{fig:alpha_functionals}, 
$\kappa_2^\star(\Palpha\invhalf A\Palpha\invhalf)$, 
$\BregmanLogDet(A, \Palpha)$, $\ln\Big(\KCond{\Palpha\inv A}\Big)$ are
plotted as a function of $\alpha$. The condition number is
flat on the interval in \cref{eq:alpha:interval}, in agreement with
\cref{thm:alpha:cond}. $\BregmanLogDet(A, \Palpha)$ and $\ln\Big(\KCond{\Palpha\inv A}\Big)$ have
a minimum at $\alpha^\star$, in agreement with \cref{thm:alpha}. Indeed,
as stated in \cref{eq:alpha}, the value of $\alpha$ that minimises Kaporin's
condition number is the \emph{average} of the remaining eigenvalues not
selected by the Bregman truncation.

\iffalse
Finally, we present a result on the influence of $\alpha$ on PCG.

\begin{theorem}\label{thm:iteration_alpha}
Let $A\in\Snpp$,  $\Palpha = \Calpha\Calpha\Hermitian$ be given by \cref{eq:Bregman_preconditioner:alpha}, 
and consider the symmetrically preconditioned system
\begin{subequations}
\begin{align}
& \Calpha\inv A \Calpha\invHermitian z = \Calpha\inv b,\label{eq:preconditioned_system:Palpha}\\
& \Calpha\Hermitian x = z.
\end{align}
\end{subequations}
The $k^\textnormal{th}$ iterate of CG applied to
\cref{eq:preconditioned_system:Palpha}, denoted by $z_k$, is 
independent of the choice of $\alpha$.
\end{theorem}
\begin{proof}
See \cref{app:iteration_proof}
\end{proof}

\Cref{thm:iteration_alpha} tells us that while computing $\alpha$ allows
us to recover Kaporin's condition number, and therefore the associated bounds on PCG convergence, it does not -- in theory -- influence PCG 
The approaches in \cite{bock2023preconditioner,bock2023new} can therefore
be viewed as similar to $\KCond{\cdot}$-minimising strategies.
In practice (i.e. finite precision
arithmetic), however, $\alpha$ likely affects the iteration history of PCG.
\fi

\section{A general result on trace scaling}\label{sec:general_trace_scaling}

We conclude with a result similar to, but more general than \cref{thm:alpha}. 
We recall that PCG converges in a single iteration when the preconditioned matrix 
is some scalar multiple of the identity, say,
\[
P\inv A = c I,
\]
for $c>0$. This motivates the study of preconditioners scaled by some constant. In particular,
we can choose the constant $c$ such that
\[%begin{equation}\label{eq:trace_scaling}
\trace{(c P)\inv A} = c\inv \trace{P\inv A} = n \Leftrightarrow c = \frac{\trace{P\inv A}}{n},
\]%end{equation}
from which we deduce
\begin{equation}\label{eq:BregmanEqualsKaporin}
\BregmanLogDet(A, cP) = \ln \KCond{(cP)\inv A}
\end{equation}
as shown in \cref{thm:Bregman_equals_lnKaporin}.

This leads to the following observation, namely that minimisation of 
$\BregmanLogDet(A, P)$ is equivalent to minimisation of $K(P\inv A)$ when
considering preconditioners for which the divergence cannot be reduced by
scaling $P$ by a positive constant.

\begin{proposition}\label{prop:DivEqualsKaporin}
Let $A\in\Snpp$ and $C\subseteq\Snpp$, and consider the 
two following optimisation problems:

\begin{equation}\label{eq:DivEqualsKaporin:Div}
S_1 = \arg\min_{P\in C} \quad \BregmanLogDet(A, P).
\end{equation}
\begin{equation}\label{eq:DivEqualsKaporin:Kaporin}
S_2 = \arg\min_{P\in C} \quad \ln \KCond{P\inv A}.
\end{equation}

If $C$ is a cone, then the following statements hold.
\begin{enumerate}[(i)]
    \item $\min_{P\in C} \BregmanLogDet(A, P) = \min_{P\in C} \ln \KCond{P\inv A},$
\item $S_1 \subset S_2$,
\item $\cone(S_1) = S_2$, where $\cone(S_1)$ denotes the conic hull of $S_1$.
\end{enumerate}
\end{proposition}
\begin{proof}
Suppose $P^\star$ is a minimiser of \cref{eq:DivEqualsKaporin:Kaporin}.
Since Kaporin's condition number is invariant under any positive scaling of
the preconditioner, we can choose $P^\star$ such that $\trace{(P^\star)\inv A} = n$.
By \cref{thm:Bregman_equals_lnKaporin}, $\BregmanLogDet(A, P^\star) = \KCond{(P^\star)\inv A}$
which proves (i). To see (ii), note that any minimiser $P$ of \cref{eq:DivEqualsKaporin:Div}
must satisfy $\trace{P\inv A} = n$. Indeed, suppose $P$ is a minimiser of 
\cref{eq:DivEqualsKaporin:Div} but $\trace{P\inv A} \neq n$. Define
\[
P_c = c P, \quad P \in S_1,
\]
where $c = \frac{\trace{P\inv A}}n$ such that $\trace{P_c\inv A} = n$. Then,
\[
\BregmanLogDet(A, P) > \BregmanLogDet(A, P_c)
\]
since
\[
\BregmanLogDet(A, P) - \BregmanLogDet(A, P_c) = n(c - 1 - \ln(c)) > 0
\]
using $x \geq 1 + \ln x$, $x>0$. As a result, $c=1$ if $P\in S_1$, since $C$ is a cone.
Further, $\BregmanLogDet(A, P) = \ln \KCond{P\inv A}$. $S_2$ contains all positive 
scalings of $P^\star$, which includes $P \in S_1$. Finally, (iii) holds by the scaling
invariance of Kaporin's condition number.
\end{proof}

We also highlight that $\KCond{P_c\inv A} = \KCond{P_{\alpha^\star}\inv A}$
for any $c>0$. In other words, $\ln\KCond{P_c\inv A}$ is a flat line on
\cref{fig:alpha_functionals} whereas $\ln\KCond{P_{\alpha}\inv A}$ 
has curvature away from $\alpha$ (recall $P_{\alpha}$ is given in
\cref{eq:Bregman_preconditioner:alpha}). This is of course obvious from 
\cref{eq:Kaporin:properties:cscale}, but there are preconditioners prevalent
in the literature, such as the randomised Nystr\"om preconditioner 
\cite{frangella2023randomized}, where a scaling similar to that controlled
by $\alpha$ is used to control the condition number of the preconditioned
matrix. Partial scaling of the
preconditioner influences the iterates of PCG (as opposed to scaling the
entire preconditioner, e.g. $P\mapsto cP$, $c>0$), so it would be interesting
and useful to analyse the consequence of a given choice beyond controlling
the condition number.

\section{Summary and future directions}\label{sec:summary}

We explored a lesser-known condition number, the Kaporin
condition number, and established its connection with the Bregman 
log determinant divergence. The link between Kaporin's condition 
number and this divergence lies in the trace scaling of the preconditioned
matrix $M$ ensuring $\frac 1n\trace{M} = 1$. This led to the error
bounds presented in \cref{eq:Kcond_logdet}. In \cref{sec:precond_low-rank}
we explore this connection in the context of the truncation described in 
\cite{bock2023new}. \Cref{sec:alphacond} described
this scaling via $\alpha$ (by scaling the orthogonal complement of the 
truncation) or $c$ (the entire preconditioner). We also highlight a
relationship with the condition number. \Cref{sec:general_trace_scaling}
contained a result concerning the equivalence of using the divergence and 
logarithm of the Kaporin condition number as an objective function.\\

% Future directions
We comment on some possible future directions.
An obvious extension would be to compute approximate factor $Q$
of a matrix $A$ using an FSAI approach \cite{janna2015fsaipack,schafer2021sparse},
and then compensate  using a low-rank term for the error $E = A - QQ\transp$ 
using the low-rank truncation provided in \cite{bock2023new}. This could lead
to a more principled approach than in \cite{bock2023new}, where $Q$ is obtained
as an incomplete Cholesky factor of $A$.\\

Randomised linear algebra is a family of techniques in 
linear algebra that leverages randomness to approximate large-scale problems that may 
otherwise be computationally infeasible \cite{halko2011finding,martinsson2020randomized}.
Such methods could prove to be useful in estimating a randomised version of the 
log determinant divergence, and therefore providing operational approximate convergence
bounds.\\

The recent paper by Carson et al. \cite{carson2024towards} contains a wealth of
insights on the behaviour of CG. The authors point out that not only is clustering 
of the eigenvalues of a matrix important for CG, but so is their position.
As described in \cite[Section 2.2]{carson2024towards}, the position of the clusters
can have an enormous effect on the residual error in energy norm. The scaling of the
preconditioned matrix via $\alpha$ given in \cref{sec:alphacond} can be thought of
as controlling the position of the "unpreconditioned" part of the matrix $A$, and it
would be interesting to explore theoretically and quantitatively the impact of $\alpha$.\\

The Bregman log determinant divergence has been studied in detail in the field
of information geometry \cite{amari2016information}.
When $N$ is a smooth submanifold of a manifold $M$, the point
dual geodesic projection $P_N^\star$ of $P$ onto $N$ is given by:
\[
P_N^\star = \arg\min_{R\in N} \Divergence (P, R),
\]
%and in the case of flat (sub)manifolds (e.g. $\Snpp$), uniqueness of such projections holds. 
In certain cases, uniqueness of such projection holds, see \cite{amari2016information}
and \cite[Section 3.7]{nielsen2020elementary} for more details. Candidate preconditioners
described by some adequate submanifold $N\subset\Snpp$ could therefore
be interesting, and insights or algorithms from information geometry 
may be relevant to preconditioning or other problems in numerical linear algebra.

\section*{Acknowledgements}
This work was supported by the Novo Nordisk Foundation under grant number NNF20OC0061894.

\appendix
\section{Proof of \cref{cor:error}}\label{app:error_proof}

\begin{proof}
\begin{align*}
\frac1n \BregmanLogDet(A,P)
%& = \frac1n \Big(\big(\sumin 1 + \lambda_i(\tilde E) - \ln(1 + \lambda_i(\tilde E))\big) - n\Big)\\
%& = \frac1n \Big(\sumin 1 + \lambda_i(\tilde E) - \ln(1 + \lambda_i(\tilde E)) - 1\Big)\\
& = \frac1n \Big(\sumin \lambda_i(\tilde E) - \ln(1 + \lambda_i(\tilde E))\Big)\\
& = \Big(\frac1n \sumin \lambda_i(\tilde E)\Big) - \frac1n \sumin \ln(1 + \lambda_i(\tilde E))\\
& \approx \ln( 1 + \frac1n \sumin \lambda_i(\tilde E)) - \frac1n \ln\det(I + \tilde E)\\
& = \ln(\frac1n \sumin 1 +  \lambda_i(\tilde E)) - \frac1n \ln\det(I + \tilde E)\\
& = \ln(\frac1n \trace{I + \tilde E}) - \frac1n \ln\det(I + \tilde E).
\end{align*}
Therefore
\begin{align*}
\BregmanLogDet(A,P) 
& \approx n \ln(\frac1n \trace{I + \tilde E}) - \ln\det(I + \tilde E)\\
& = n \ln(\frac1n \trace{I + \tilde E}) - \ln\det(I + \tilde E)\\
& = \ln (K(M)^n),
\end{align*}
and
\[
\ln (K(M)^n) = \BregmanLogDet(A,P) + O(n^2 \trace{\tilde E}).
\]
\end{proof}

\bibliographystyle{siam}
\bibliography{main}
%\printbibliography
\end{document}